\documentclass{article}
\usepackage{amsmath,amssymb,amsfonts,amsthm}
\usepackage{graphicx}
\usepackage{layout}
\usepackage[latin1]{inputenc}
\usepackage[french,english]{babel}
\usepackage[french]{babel}
\usepackage{cancel}

\newtheorem{theorem}{Theorem}

\newtheorem{corollary}{Corollary}
\newtheorem{proposition}{Proposition}
\newtheorem{lemma}{Lemma}

\newtheorem{notation}{Notation}

\newtheorem{example}{Example}
\newtheorem{remark}{Remark}    

\usepackage{fancyhdr}
\usepackage{setspace}

\renewcommand{\headrulewidth}{0,5pt}
\newcommand{\vertiii}[1]{{\left\vert\kern-0.25ex\left\vert\kern-0.25ex\left\vert #1 
    \right\vert\kern-0.25ex\right\vert\kern-0.25ex\right\vert}} %Norme trois lignes
\title{On characterizations of a some classes of Schauder frames in Banach spaces}
\author{{\footnotesize Rafik Karkri, Samir Kabbaj and Hamad Sidi Lafdal}  }

\date{\today}
\begin{document}
\maketitle %pour afficher le titre et les auteurs ...

\begin{abstract}
In this paper, we prove the following results. There exists a  Banach space
without  basis which has a Schauder frame.
There  exists an universal Banach space $B$  (resp. $\tilde{B}$) 
with a  basis (resp. an unconditional  basis)  such that, a  Banach  $X$
has a Schauder frame (resp.  an unconditional  Schauder frame )
if and only if $X$ is isomorphic 
to a complemented subspace of $B$ (resp. $\tilde{B}$).
For a weakly sequentially complete Banach space, a Schauder frame is unconditional
if and only if it is besselian. 
 A separable Banach space $X$ has a Schauder frame if and only if
it has the bounded approximation property. 
Consequenty, The Banach space $\mathcal{L}(\mathcal{H},\mathcal{H})$ of all 
bounded linear operators on a Hilbert space $\mathcal{H}$  has no Schauder frame.
Also, if $X$ and $Y$ are Banach spaces  with Schauder frames then, the
 Banach space $  X\widehat{\otimes}_{\pi}Y$
 (the projective tensor product of $X$ and $Y$) has a Schauder frame.
  From the Faber$-$Schauder system  we  construct a Schauder frame
for the Banach space $C[0,1]$ 
(the Banach space of continuous functions on the closed
interval $ [0,1]$)  which is not a  Schauder basis of $C[0,1]$.
Finally, we give a positive answer to some open problems related 
 to the Schauder bases (In the Schauder frames setting).
  
\end{abstract}
{\footnotesize $\textbf{MSC2020-Mathematics Subject Classification System}$. 
46B04, 46B10, 46B15, 46B25, 46B45.}\\
{\footnotesize $\textbf{Keywords and phrases:}
$ Schauder frames, Besselian Schauder frames, Unconditional Schauder frames,
 Weakly sequentially complete Banach spaces, Universal Banach space,
 Bounded approximation property.}

\renewcommand{\headrulewidth}{0pt}  % epaisseur  de la ligne en bas et pied du page
%\fancyfoot[C]{{\scriptsize R. Karkri,  S. Kabbaj and H. Sidi Lafdal}}
\fancyhead[C]{{\scriptsize R. Karkri,  S. Kabbaj and H. Sidi Lafdal}} %afficher en bas à droite de tous les page  

\maketitle

\section{Introduction}
In 1932 Banach and Mazur \cite[theorem 9, p.185]{S.Banach.1932} 
state that every separable Banach space is isometrically embeddable in $C[0,1]$.
The fact that $C[0,1]$ has a  Schauder basis \cite{J.Schauder},\cite[example 4.1.11, p.352]{meg}
 led Banach to
pose the question whether every separable Banach space has a Schauder basis 
\cite[p.111]{S.Banach.1932}. We recall that the sequence
$\left\lbrace x_{n} \right\rbrace _{n\in \mathbb{N}^{*}}$ in a Banach space 
$X$ is called a Schauder basis of $X$ if for every $x\in X$ there is
a unique sequence of scalars $\left\lbrace \alpha_{n} \right\rbrace _{n\in \mathbb{N}^{*}}$
so that 
\[x=\underset{n=1}{\overset{+\infty }{\sum }}\alpha_{n}x_{n}\]
This problem (known as the basis problem) remained open for a long time and was solved in 1973 in the
negative by P. Enflo \cite{P.Enflo.1973}.
In 1969   A. Pe\l czy\'{n}ski \cite{pel2} constructed a complementably universal Banach space for
the class of Banach spaces with a Schauder basis (theorem \ref{uni.basis.1969.C1}).
In 1987 S. J. Szarek \cite{S.J.Szarek} showed that there is a complemented subspace
 without Schauder basis of a space with a Schauder basis and answered therefore
to a problem of fifty years old.
That is the property of having a Schauder basis is not inherited by
complemented subspaces. 

The purpose of this paper is to see what happens if we replace Schauder bases by 
Schauder frames in the above results and questions. 
In the following, we give a brief history of frame theory.

In $1946,$ Gabor \cite{gab} formulated a fundamental approach 
to signal decomposition in terms of
elementary signals. In $1952,$ Duffin and Schaeffer  
\cite{duf} abstracted the fundamental notion of Gabor for
	  introduced  the notion of frames for Hilbert spaces as a generalization 
	  of bases.
  A frame for a Hilbert space  
 $\left( \mathcal{H},\left \langle .,.\right \rangle \right) $ is a sequence
 $\left\lbrace h_{n}\right\rbrace _{n\in \mathbb{N}^{* }}$  of elements of 
 $\mathcal{H}$ for which  there exists  constants $A,B>0$ such that: 
\begin{equation}
\label{std.frame}
A\left \Vert h\right \Vert _{\mathcal{H}}^{2}\leq \underset{n=1}{\overset{+\infty }{\sum }}
\left \vert \left \langle h,h_{n}\right \rangle \right \vert ^{2}
\leq B\left \Vert h\right \Vert _{\mathcal{H}}^{2},\; \;h\in \mathcal{H}
\end{equation}
The numbers $A$ and $B$ are called bounds for the frame.
If the upper inequality in \eqref{std.frame} holds,
$\left\lbrace h_{n}\right\rbrace _{n\in \mathbb{N}^{*}}$ is said to be 
a Bessel sequence with bound $B$. This condition allows each element in  
$\mathcal{H}$ to be 
written as a linear combination of the elements of 
$\left\lbrace   h_{n}\right\rbrace  _{n\in \mathbb{N}^{*}}$, but  linear
independence between the elements of 
 $\left\lbrace   h_{n}\right\rbrace  _{n\in \mathbb{N}^{*}}$ 
 is not required. 
 A sequence $\left\lbrace   h_{n}\right\rbrace  _{n\in \mathbb{N}^{*}}$
 in $\mathcal{H}$  is bounded if 
 $ 0<inf_{n\in \mathbb{N}^{*}}\Vert h_{n}\Vert_{\mathcal{H}}
 \leq sup_{n\in \mathbb{N}^{*}}\Vert h_{n}\Vert_{\mathcal{H}}<\infty
  $.
A Riesz basis for $\mathcal{H}$ is a family of the form
$\left\lbrace Ue_{n}\right\rbrace_{n\in \mathbb{N}^{* }}$,
where  $\left\lbrace e_{n}\right\rbrace_{n\in \mathbb{N}^{*}}$
is an orthonormal basis for $\mathcal{H}$ and
 $U : \mathcal{H}\longrightarrow \mathcal{H}$
 is a bounded bijective operator \cite[p.86]{O.Christensen.2016}. 
 A sequence $\left\lbrace   h_{n}\right\rbrace  _{n\in \mathbb{N}^{*}}$ 
 which is a Riesz basis for its closed linear span in $\mathcal{H}$ 
 is called a Riesz basic sequence in $\mathcal{H}$.
In 2013 \cite{M.S.N.S.2015.Fei.Conjec}  Marcus, Spielman,
and Srivastava confirmed the Feichtinger conjecture
(Every bounded frame can be written as a finite union of Riesz basic sequences).

 A frames for a Banach spaces was 
 introduced in $1991$ by Gr\"{o}chenig \cite{gro}
 to extend the definition of frames to that of general separable Banach spaces. 
 In light of the works of  Cassaza, Han and Larson  \cite
{cas} $(1999)$,    Han and Larson \cite{lar2} $(2000)$.
 Cassaza \cite{cas.2008} $(2008)$  introduced the notion of
Schauder frame of a given  Banach space. 
A Schauder frame of a given  Banach spaces $X$ is  a sequence  
$\left\lbrace (x_{n},f_{n})\right\rbrace _{n\in \mathbb{N}^{*}}\subset X\times X^{*}$
 such that, for each $x$ in $X$ we have:
 \[ x=\underset{n=1}{\overset{+\infty}{\sum }}
 f_{n}\left(x\right) x_{n}.\]

In 2021 Karkri and Zoubeir
  \cite{karkri.zoubeir.2021} introduced the notion of
besselian Schauder frames of Banach spaces and obtained in another way
the existence of an universal space  for the class of separable Banach  spaces
with Schauder frames (corollary \ref{characterization.conditional.1}). Also, proved 
(corollary \ref{universal.2})
for the class of separable weakly sequentially complete Banach  spaces 
 with  besselian Schauder frames.

\section{Preliminaries}
Let $\left( X,\left \Vert .\right \Vert _{X}\right) $ 
be a Banach space on  $\mathbb{K}\in \left \{ \mathbb{R},\mathbb{C}\right \} $ 
 and $X^{*}$ it's
topological dual.

\begin{enumerate}
\item We denote by $\mathbb{B}_{X}$ the closed unit ball of $X:$
\begin{equation*}
\mathbb{B}_{X}:=\left \{ x\in X:\left \Vert x\right \Vert _{X}\leq 1\right \}
\end{equation*}
\item
$C[0,1]$ is the Banach space of continuous functions on the closed
interval $ [0,1]$.
\item 
Let $1< p < \infty $. We denote by $L_{p} [0,1]   $ the Banach space of equivalence
 classes of Lebesgue integrable functions 
 $f: [0,1]\longrightarrow \mathbb{K}$, with the norm
 \[ \Vert f\Vert_{L_{p} [0,1]   }=\left(\int_{0}^{1}\vert f \vert^{p} d\mu\right)^{1/p}\]
\item
\cite[p.57]{W.B.Johnson.J.Lindenstrauss.2001} Let $1< p < \infty$ and $\lambda\geq 1$.
A Banach space $X$ is called an $\mathcal{L}_{p,\lambda}$
 space if for every finite dimensional subspace $E$ of $X$
  there is a further subspace $ F\supset E$ with $d(F,l_{p}^{n})\leq \lambda$
  where $n=dim F$. A space is called an $\mathcal{L}_{p}$ space
   if it is an $\mathcal{L}_{p,\lambda}$  space for some $ \lambda <\infty $.
   Where $l_{p}^{n}$ is  the $\mathbb{K}$-vector space of
sequences $\lambda :=\left\lbrace \lambda _{k}\right\rbrace _{1\leq k\leq n}$ 
of scalars, endowed with the norm
$\left \Vert \lambda \right \Vert 
=\left (\sum_{k=1}^{n}\left \vert \lambda _{k}\right \vert^{p}\right )^{1/p} $.
\item
\cite[p.5]{W.B.Johnson.J.Lindenstrauss.2001} 
Let $ \left \lbrace X_{n}\right\rbrace _{n\in \mathbb{N}^{*}}  $
be a sequence of Banach spaces and $1< p < \infty $. 
We denote by $\left (\sum_{n=1}^{\infty}X_{n}\right )_{p}$ 
the Banach space of all sequences 
$ \left\lbrace x_{n}\right\rbrace _{n\in \mathbb{N}^{*}}  $ with 
$x_{n}\in X_{n}$ and
$\left\Vert \left\lbrace x_{n}\right\rbrace _{n\in \mathbb{N}^{*}} \right \Vert
=\left (\sum_{n=1}^{\infty}\left\Vert x_{n}\right \Vert^{p}\right )^{1/p}$.
  
\item For each $\left( i,j\right) \in \mathbb{N}^{2},$ we denote by $\delta _{i,j}$
 the so-called Kroneker symbole defined by 
\begin{equation*}
\delta _{i,j}:=\left \{ 
\begin{array}{cc}
1& if\; i=j \\ 
0& else
\end{array}
\right.
\end{equation*}
\item
 We denote by 
$\mathfrak{S}_{\mathbb{N}^{*}}$ the set of all permutations of $\mathbb{N}^{* }$.

\item
We denote by  $\mathcal{D}_{\mathbb{N}^{*}}$
 the set of all the nonempty finite subsets of $\mathbb{N}^{*}$.

\item For each $L\in \mathcal{D}_{\mathbb{N}^{*}}$,
 we denote by $\min\left( L\right) $ the minimal element of $L$.

\item The Banach space $X$ is said to be weakly sequentially complete 
(see \cite[definition 2.5.23, p.218]{meg} and \cite[p.37-38]{kal})  if for each
sequence $\lbrace x_{n}\rbrace _{n\in \mathbb{N}^{* }}$ of $X$ 
such that $\underset{n\rightarrow +\infty }{\lim }f(x_{n})$ 
exists for every $f\in X^{*}$, there exists $x\in X$ such that
 $$\underset{n\rightarrow+\infty }{\lim }f(x_{n})=f(x)\;\;\;(f\in X^{*})$$
\item A sequence $\mathcal{F}:=\left\lbrace \left( x_{n},f_{n}\right)\right\rbrace _{n\in \mathbb{N}^{*}}\subset X\times X^{*}$ is called a paire of $X$.

\item The paire $\mathcal{F}$ is called a Schauder frame (resp.
unconditional Schauder frame ) of $X$ if for all $x\in X$, the series 
$\sum f_{n}\left( x\right) x_{n}$ is convergent (resp. unconditionally
convergent) in $X$ to $x$.

\item The paire $\mathcal{F}$ is said to be a besselian paire of $X$ if
there exists a constant $C>0$ such that : 
\begin{equation*}
\underset{n=1}{\overset{+\infty }{\sum }}
\left \vert f_{n}\left(x\right) \right \vert \left \vert f\left( x_{n}\right) \right \vert
\leq C\left \Vert x\right \Vert _{X}\left \Vert f\right \Vert
_{X^{* }}
\end{equation*}
for each $x\in X$ and $f\in X^{*}.$

\item The paire $\mathcal{F}$ is said to be a besselian Schauder frame of $X$ 
if it is both a besselian paire and a Schauder frame of $X$ .
\end{enumerate}
\begin{remark}
For a besselian paire $\mathcal{F}$ of $X$, the quantity 
\begin{equation*}
\mathcal{L}_{\mathcal{F}}:=\underset{\left( x,f\right) \in \mathbb{B}_{X}\times \mathbb{B}_{X^{*}}}{\sup }\left( \underset{n=1}{\overset{%
+\infty }{\sum }}\left \vert f_{n}\left( x\right) \right \vert \left
\vert f\left( x_{n}\right) \right \vert \right)
\end{equation*}
is finite and for each $ x\in X$ and $f\in X^{* }$,
 the following inequality holds
\begin{equation*}
\underset{n=1}{\overset{+\infty }{\sum }}
\left \vert f_{n}\left(x\right) \right \vert 
\left \vert f\left( x_{n}\right) \right \vert
\leq \mathcal{L}_{\mathcal{F}}\left \Vert x\right\Vert _{X}
\left \Vert f\right \Vert_{X^{*}}
\end{equation*}  
\end{remark}
\begin{remark}
The Banach spaces $L_{1}[0,1]$ has a  Schauder frame but 
has no besselian Schauder frame \cite[corollary 5.3, p.820]{karkri.zoubeir.2021}.
\end{remark}
\begin{remark}
Let $X$ be a separable Banach space. there exists  a besselian paire of $X$ 
which is not  Schauder frame of $X$.
\end{remark}
\textbf{Proof.}
Let $\left\lbrace \left( x_{n},f_{n}\right)\right\rbrace _{n\in \mathbb{N}^{*}}
$ be  a paire of $X$ 
such that
$x_{n}\neq 0$ and $f_{n}\neq 0$ for each $n\in\mathbb{N}^{*}$ 
and $\left\lbrace x_{n} \right\rbrace_{n\in\mathbb{N}^{*}}$
 is not complete in $X$.
We put for each $n\in\mathbb{N}^{*}$,
  $y_{n}=x_{n}/ \sqrt{2}^{n}\left \Vert  x_{n} \right \Vert $ and
$g_{n}=f_{n}/ \sqrt{2}^{n}\left \Vert  f_{n} \right \Vert $.
Then for each $x\in X$ and $f\in X^{*}$ we have
\begin{align*}
\underset{n=1}{\overset{+\infty }{\sum }}
\left \vert g_{n}\left(x\right) \right \vert
 \left \vert f\left( y_{n}\right) \right \vert
 &\leq \underset{n=1}{\overset{+\infty }{\sum }}
 \left \Vert g_{n}\right \Vert_{X^{*}} \left \Vert x  \right \Vert_{X} 
 \left \Vert f \right \Vert_{X^{*}} \left \Vert y_{n} \right \Vert_{X}\\
 &\leq \left \Vert x \right\Vert_{X} \left \Vert f \right \Vert_{X^{*}}
 \underset{n=1}{\overset{+\infty }{\sum }} 1/2^{n}\\
 &\leq \left \Vert x \right \Vert_{X} \left \Vert f \right \Vert_{X^{*}}
\end{align*}
Consequently, 
$\left\lbrace \left( y_{n},g_{n}\right)\right\rbrace _{n\in \mathbb{N}^{*}}$ 
is  a besselian paire
but is  not a Schauder frame of $X$.
\qed

For all the material on Banach spaces, one can refer to 
\cite{meg, kal, lin01, lin02, woj}. In the sequel 
$\left(X,\left\Vert \cdot \right \Vert _{X}\right) $ 
is a given separable Banach space and
 $\mathcal{F}:=\left\lbrace \left( x_{n},f_{n}\right) \right\rbrace
_{n\in \mathbb{N}^{*}}$ is a fixed paire of $X$ 
such that $x_{n}\neq 0$ $(n=1,2,...)$.
\begin{example}
\label{Faber.frame}
There exists a sequence of piecewise linear
functions $\left\lbrace \varphi_{n} \right\rbrace _{n\in \mathbb{N}^{*}}$ in $C[0,1]$
and a sequence $\left\lbrace A_{n} \right\rbrace _{n\in \mathbb{N}^{*}}$ 
in $C[0,1]^{*}$ such that
$\left\lbrace \left( A_{n},\varphi_{n}\right) \right\rbrace
_{n\in \mathbb{N}^{*}}$ is a Schauder frame of $C[0,1]$ but is not a Schauder basis.
\end{example}
\begin{notation}
Let  $m\in \mathbb{N}$ and $0\leq k \leq 2^{m}-1$. We put \\
$J_{m}=\lbrace 0,1/2^{m},2/2^{m},3/2^{m},...,(2^{m}-1)/2^{m},1\rbrace $ and
$I_{m,k}=[ k/2^{m},(k+1)/2^{m} ]$.
\end{notation}
\begin{lemma}
There exists a sequence of piecewise linear
functions $\left\lbrace \varphi_{n} \right\rbrace _{n\in \mathbb{N}^{*}}$ in $C[0,1]$
and a sequence $\left\lbrace A_{n} \right\rbrace _{n\in \mathbb{N}^{*}}$ in $C[0,1]^{*}$
such that
\begin{align}
\label{series.rational}
f(p)=S_{2^{m}}(f)(p)=S_{2^{m}+k}(f)(p)
\;\;\;\;
 \left (m\in \mathbb{N}^{*}, 1\leq k\leq 2^{m}-1,  p\in J_{m-1}, f\in C[0,1]\right )
\end{align}
where
\begin{enumerate}
\item
$S_{2^{m}}(f)=A_{1}(f)\varphi_{1}+\underset{n=0}{\overset{m-1}{\sum }}
\underset{i=1}{\overset{2^{n}}{\sum }}
A_{2^{n}+i}(f)\varphi_{2^{n}+i}\;\;\;\; (f\in C[0,1], m\geq 1)$.
\item
For each $m\in \mathbb{N}^{*}$ and $f\in C[0,1]$ the function 
$ S_{2^{m}}(f)$ is linear in the intervals $I_{m-1,k}$ $(0\leq k\leq 2^{m-1}-1)$.
\item
$S_{2^{m}+k}(f)=S_{2^{m}}(f)
+\underset{j=1}{\overset{k}{\sum }}
A_{2^{n}+j}(f)\varphi_{2^{n}+j}
\;\;\;\; (f\in C[0,1], m\geq 1, 1\leq k\leq 2^{m}-1)$
\item
For each $m\in \mathbb{N}^{*}$, $1\leq k\leq 2^{m}-1$ and $f\in C[0,1]$ the function 
$ S_{2^{m}+k}(f)$ is linear in the intervals $I_{m,t}$ $(0\leq t\leq 2^{m}-1)$.
\end{enumerate}
\end{lemma}
\textbf{Proof.}
Let $\left( \lambda_{n} \right) _{n\in \mathbb{N}^{*}}$ be a sequence of scalars
such that $0\leq \lambda_{n} \leq 1$ $(n\in \mathbb{N}^{*})$.
\begin{itemize}
\item For $m=1$. 
 We put for each $x\in [0,1]$ and $f\in C[0,1]$: $\varphi_{1}(x)=x$, $\varphi_{2}(x)=1-x$ 
$A_{1}(f)=f(1)$ and $A_{2}(f)=f(0)$. Then
\[f(p)=S_{2^{1}}(f)(p)\;\;\;\; (p\in J_{0}, f\in C[0,1])\]
\item For $m=2$.  We define $\varphi_{2^{1}+1} $ and  
$A_{2^{1}+1}$ such that
\begin{itemize}
\item
$\varphi_{2^{1}+1}(1/2)=1 $, $\varphi_{2^{1}+1}(p)=0 $ $(p\in J_{0}) $ 
and $\varphi_{2^{1}+1}$ is linear in the intervals
$I_{1,0}$ and $I_{1,1}$. 
\item
 $A_{2^{1}+1}(f)=\lambda_{2^{1}+1} f(0)+(1-\lambda_{2^{1}+1})f(1)
 -S_{2^{1}}(f)(1/2)$ for each $f$ in $ C[0,1]$.
\end{itemize}
We define $\varphi_{2^{1}+2} $ and  
$A_{2^{1}+2}$ such that
\begin{itemize}
\item
$\varphi_{2^{1}+2}=\varphi_{2^{1}+1} $. 
\item
 $A_{2^{1}+2}(f)=f(1/2)-S_{2^{1}+1}(f)(1/2)$ for each $f$ in $ C[0,1]$.
\end{itemize}
then
\[f(p)=S_{2^{2}}(f)(p)\;\;\;\; (p\in J_{1}, f\in C[0,1])\]
 \item
Assume that the equality \eqref{series.rational} hold for $m\geq 2$.
We define the functions 
$\left\lbrace \varphi_{2^{m}+k} \right\rbrace _{1\leq k \leq 2^{m}}$ 
and
$\left\lbrace A_{2^{m}+k} \right\rbrace _{1\leq k \leq 2^{m}}$
such that:
\begin{itemize}
\item
$\varphi_{2^{m}+k}(p)=0$ $\left (1\leq k \leq 2^{m-1}, p\in J_{m-1}\right )$. 
\item
$\varphi_{2^{m}+k}\left (\frac{2k-1}{2^{m}}\right )=1$ $\left (1\leq k \leq 2^{m-1}\right )$. 
\item
For each $1\leq k \leq 2^{m-1}$ the function $\varphi_{2^{m}+k}$ is linear in the
intervals
$I_{m,0}$,..., $I_{m,2^{m}-1}$.
\item
$\varphi_{2^{m}+2^{m-1}+k}=\varphi_{2^{m}+k}$ $\left (1\leq k \leq 2^{m-1}\right )$. 
\item
For each $1\leq k \leq 2^{m-1}$ and $f$ in $ C[0,1]$ we put
 $$A_{2^{m}+k}(f)
 =\lambda_{2^{m}+k} f\left(\frac{k-1}{2^{m-1}}\right)
 +(1-\lambda_{2^{m}+k})f\left(\frac{k}{2^{m-1}}\right)
 -S_{2^{m}+k-1}(f)\left(\frac{2k-1}{2^{m}}\right)$$
 $$A_{2^{m}+2^{m-1}+k}(f)
 =f\left(\frac{2k-1}{2^{m}}\right)
 -S_{2^{m}+2^{m-1}+k-1}(f)\left(\frac{2k-1}{2^{m}}\right)$$
\end{itemize}
then the equality \eqref{series.rational} hold for $m+1$.
\end{itemize}
\qed\\
\textbf{Proof of example \eqref{Faber.frame}.} Let $f\in C[0,1]$ and $\varepsilon>0$. Since $f$
 is uniformly continuous on $[0,1]$, there exists $m_{0}\in \mathbb{N}^{*}$ 
 such that $\left \vert f(x)-f(y)\right \vert < \varepsilon/2$ whenever 
 $x,y \in [0,1]$, $\left \vert x-y\right \vert < 1/2^{m_{0}-1}$. Now,
 let $x\in [0,1]$ be arbitrary. Let  $m\geq  m_{0}$ and $1\leq k\leq 2^{m}-1$. 
 There exists $1\leq i\leq 2^{m-1}$  
 such that $\frac{i-1}{2^{m-1}}\leq x\leq \frac{i}{2^{m-1}}$.
 We may assume, without loss of generality, that 
 $\frac{2i-1}{2^{m}}\leq x\leq \frac{i}{2^{m-1}}$. 
 Since $S_{2^{m}+k}(f)$ is linear in the interval
 $\left [\frac{2i-1}{2^{m}}, \frac{i}{2^{m-1}}\right ]$, there exists
 $0\leq \beta_{2^{m}+k}\leq 1$
  and  
  $0\leq \gamma_{2^{m}+k}\leq 1$
  such that
  \begin{align*}
  S_{2^{m}+k}(f)(x)=&\beta_{2^{m}+k}S_{2^{m}+k}(f)\left (\frac{2i-1}{2^{m}}\right )
 +(1-\beta_{2^{m}+k})S_{2^{m}+k}(f)\left (\frac{i}{2^{m-1}}\right )\\
 =&\beta_{2^{m}+k}S_{2^{m}+k}(f)\left (\frac{2i-1}{2^{m}}\right )
 +(1-\beta_{2^{m}+k})f\left (\frac{i}{2^{m-1}}\right ) 
  \end{align*}
 and
  \[S_{2^{m}+k}(f)\left (\frac{2i-1}{2^{m}}\right )
  =\gamma_{2^{m}+k}f\left (\frac{i-1}{2^{m-1}}\right )
 +(1-\gamma_{2^{m}+k})f\left (\frac{i}{2^{m-1}}\right )\]
then,
\begin{align*}
\left \vert f(x)-S_{2^{m}+k}(f)(x)\right \vert 
\leq &\left \vert f(x)-f\left (\frac{i}{2^{m-1}}\right ) \right \vert 
+\left \vert f\left (\frac{i}{2^{m-1}}\right ) 
-S_{2^{m}+k}(f)(x) \right \vert\\
\leq &\frac{\varepsilon}{2}
+\beta_{2^{m}+k}\gamma_{2^{m}+k}\left \vert f\left (\frac{i-1}{2^{m-1}}\right ) 
-f\left (\frac{i}{2^{m-1}}\right )  \right \vert\\
\leq &\frac{\varepsilon}{2}
+\frac{\varepsilon}{2}\beta_{2^{m}+k}\gamma_{2^{m}+k}\\
\leq &\varepsilon
\end{align*}
Consequently, 
$$\underset{n\longrightarrow \infty}{lim}
\left \Vert f-\underset{j=1}{\overset{n}{\sum }}
A_{j}(f)\varphi_{j}  \right \Vert_{C[0,1]}
=\underset{n\longrightarrow \infty}{lim}
\underset{x\in[0,1]}{sup}
\left \vert f(x)-\underset{j=1}{\overset{n}{\sum }}
A_{j}(f)\varphi_{j}(x)  \right \vert=0$$
\qed

\section{Schauder frames}
\begin{lemma}\cite[remark 7.1, p.189]{I.Singer.II}
Let $X$ be a Banach space and $  \left\lbrace f_{n}\right\rbrace _{n\in \mathbb{N}^{*}}\subset X^{*} $ be a
total sequence of functionals. Then $X$ is linearly isometric to the Banach space 
$$S_{X}=\left\lbrace\left\lbrace  f_{n}(x)\right\rbrace_{n\in\mathbb{N}^{*}}:
 x\in X\right\rbrace$$
(where the norm is defined by 
$ \left \Vert \left\lbrace  f_{n}(x)\right\rbrace_{n\in\mathbb{N}^{*}} \right \Vert_{S_{X}}
=\left \Vert x \right \Vert_{X}$), by the mapping
$$ \psi : x\longmapsto \left\lbrace  f_{n}^{*}(x)\right\rbrace_{n\in\mathbb{N}^{* }}$$
\end{lemma}
\begin{remark}
\label{frame.operators.spaces}
Let $\mathcal{F}$ 
be a Schauder frame of a Banach space $X$. Then
$\left\lbrace \left(\psi(x_{n}),f_{n}\circ \psi ^{-1}\right) \right\rbrace _{n\in \mathbb{N}^{* }}$
is a Schauder frame of a Banach space $S_{X}$.
\end{remark}
\begin{proposition}\cite[proposition 3.1, p.18]{I.Singer.I}
\label{B.space.A}
Let $\left\lbrace x_{n} \right\rbrace_{n\in\mathbb{N}^{* }}$
be a sequence in a Banach space $X$, 
such that $x_{n}\neq 0$ $(n=1,2,...)$ and let 
$\mathcal{A}$ be the linear space
of sequences of scalars
\[\mathcal{A}=\left\lbrace
 \left\lbrace \alpha_{n} \right\rbrace_{n\in\mathbb{N}^{*}}\subset \mathbb{K}:
\sum _{k}\alpha_{k}x_{k}\;converges \right\rbrace\]
endowed with the norm
\[ \left \Vert \left\lbrace \alpha_{n} \right\rbrace_{n\in\mathbb{N}^{*}}
 \right \Vert_{\mathcal{A}}
=\underset{n\in\mathbb{N}^{*}}{\sup }
\left \Vert \overset{n}{\underset{k=1}{\sum }} \alpha_{k}x_{k}  \right \Vert_{X}\]
Then $\mathcal{A}$ is a Banach space.
\end{proposition}
\begin{remark}
If   $\left\lbrace x_{n} \right\rbrace_{n\in\mathbb{N}^{*}}$
is a sequence in a Banach space $X$ such that $x_{n}\neq 0$ $(n=1,2,...)$, 
then the mapping 
\begin{align*}
T:\mathcal{A} &\rightarrow  X \\ 
 {\left\lbrace \alpha_{n} \right\rbrace}  &\mapsto \overset{+\infty }{\underset{n=1}{\sum }}
 \alpha_{k}x_{k}
\end{align*}
defines a bounded linear operator.
\end{remark}
\textbf{Proof.}
For each $\alpha=\left\lbrace \alpha_{n} \right\rbrace_{n\in\mathbb{N}^{*}}
\in \mathbb{B}_{\mathcal{A}}$ we have
 $\left \Vert \sum_{k=1}^{n}\alpha_{k}x_{k} \right \Vert_{X}\leq 1 \;\;\;(n\in\mathbb{N}^{*})$.
 Consequently,
\begin{align*}
\left \Vert T \right \Vert
&=\underset{\alpha\in\mathbb{B}_{\mathcal{A}}}{\sup }
\left \Vert T\alpha \right \Vert_{X}\\
 &=\underset{\alpha\in\mathbb{B}_{\mathcal{A}}}{\sup }
\left \Vert \overset{\infty}{\underset{k=1}{\sum }} \alpha_{k}x_{k}  \right \Vert_{X}\\
&\leq 1
\end{align*}
\qed
\begin{remark}
If  $\mathcal{F}$ 
is a Schauder frame of a Banach space $X $ 
 and $k\in\mathbb{N}^{*}$, then the  mapping 
\begin{align*}
T_{k}: \mathcal{A} &\rightarrow  X \\ 
 \alpha={\left\lbrace \alpha_{n} \right\rbrace}  &\mapsto \overset{k}{\underset{i=1}{\sum }}
 f_{i}(T\alpha)x_{i}
\end{align*}
defines a bounded linear operator.
\end{remark}
\textbf{Proof.}
For each $k\in\mathbb{N}^{*}$, $T_{k}$ is a finite sum of bounded linear operators
$$\alpha \longmapsto f_{i}(T\alpha)x_{i}, \;\;\;(1\leq i\leq k)$$\qed
\begin{remark}
Assume that $\mathcal{F}$ is a Schauder frame of $X$. Then
the Banach space $S_{X}$ is a subspace of $ \mathcal{A}$. 
\end{remark}
\begin{proposition}\cite[proposition 8.1, p.74]{I.Singer.I}
The unit vectors 
$e_{n}=\left\lbrace \delta_{n,k} \right\rbrace_{k\in\mathbb{N}^{* }}$ 
$(n=1,2,...)$ constitue a Schauder basis of  $\mathcal{A}$.
\end{proposition}
\begin{theorem}
\label{projection.onto}
Let $\mathcal{F}$ 
be a Schauder frame of a Banach space $X $. Then 
the infinite matrix
\[\mathcal{M}_{i,j}=f_{i}\left( x_{j}\right ),\; (i,j\in \mathbb{N}^{*})\]
defines a bounded projection $\mathcal{M}$ from $\mathcal{A}$ onto $S_{X}$. 
\end{theorem}
\textbf{Proof.}
Let 
$ \alpha=\left\lbrace \alpha_{n} \right\rbrace_{n\in\mathbb{N}^{*}}$
be a sequence in $\mathcal{A}$. For each $i\in  \mathbb{N}^{*}$ we have
\begin{align*}
f_{i}\left(T\alpha \right )&=f_{i}\left(\overset{+\infty }{\underset{k=1}{\sum }}
 \alpha_{k}x_{k}\right )\\
 &=\overset{+\infty }{\underset{k=1}{\sum }}
f_{i}\left( x_{k}\right )\alpha_{k}\\
& =\overset{+\infty }{\underset{k=1}{\sum }}
\mathcal{M}_{i,k}\alpha_{k}
\end{align*}
that is 
$\mathcal{M}\alpha=\left\lbrace f_{n}\left(T\alpha \right ) \right\rbrace_{n\in\mathbb{N}^{*}}$.
Since $T\alpha=\overset{+\infty }{\underset{k=1}{\sum }}
f_{k}\left(T\alpha \right )x_{k}$, it follows that
 $\mathcal{M}\alpha\in \mathcal{A}$.
 On the other hand, for each $i,j\in  \mathbb{N}^{*} $ we have
\begin{align*}
\left(\mathcal{M}^{2} \right )_{i,j}&=\overset{+\infty }{\underset{k=1}{\sum }}
\mathcal{M}_{i,k}\mathcal{M}_{k,j}\\
 &=\overset{+\infty }{\underset{k=1}{\sum }}
f_{i}\left( x_{k}\right )f_{k}\left( x_{j}\right ) \\
& =f_{i}\left(\overset{+\infty }{\underset{k=1}{\sum }}
 f_{k}\left( x_{j}\right )x_{k}\right ) \\
 & =f_{i}\left(x_{j}\right )\\
& =\mathcal{M}_{i,j}
\end{align*}
then, $\mathcal{M}$ is a projection in $\mathcal{A}$.
By  the
uniform boundedness principle applied  to the family of bounded linear operators
$\left\lbrace T_{n} \right\rbrace_{n\in \mathbb{N}^{*}}$ we have
 \begin{align*}
 \left \Vert \mathcal{M} \right \Vert&=\underset{\alpha\in\mathbb{B}_{\mathcal{A}}}{\sup }
 \left \Vert \mathcal{M}\alpha \right \Vert_{\mathcal{A}}\\
 &=\underset{\alpha\in\mathbb{B}_{\mathcal{A}}}{\sup }
 \left \Vert\left\lbrace f_{n}\left(T\alpha \right ) \right\rbrace_{n\in\mathbb{N}^{*}}\right \Vert_{\mathcal{A}}\\
&=\underset{\alpha\in\mathbb{B}_{\mathcal{A}}}{\sup }
 \underset{n\in\mathbb{N}^{*}}{\sup }
 \left \Vert \overset{n}{\underset{k=1}{\sum }}
f_{k}\left(T\alpha \right )x_{k} \right \Vert_{X}\\
&=\underset{\alpha\in\mathbb{B}_{\mathcal{A}}}{\sup }
 \underset{n\in\mathbb{N}^{*}}{\sup }
 \left \Vert T_{n}\alpha  \right \Vert_{X}\\
 &=\underset{n\in\mathbb{N}^{*}}{\sup }
 \underset{\alpha\in\mathbb{B}_{\mathcal{A}}}{\sup }
 \left \Vert T_{n}\alpha  \right \Vert_{X}\\
 &=\underset{n\in\mathbb{N}^{*}}{\sup }
 \left \Vert T_{n}\right \Vert <\infty
 \end{align*}
then $\mathcal{M}$ is a bounded operator. For each $i,n\in  \mathbb{N}^{*} $ we have
\begin{align*}
\left(\mathcal{M}e_{n} \right )_{i}&=\overset{+\infty }{\underset{k=1}{\sum }}
\mathcal{M}_{i,k}\delta_{n,k}\\
 &=\mathcal{M}_{i,n}\delta_{n,n}\\
& =f_{i}\left( x_{n}\right )
\end{align*}
then $\mathcal{M}e_{n}=\left\lbrace f_{i}\left(x_{n} \right ) \right\rbrace_{i\in\mathbb{N}^{*}}$.
It follows from remark \eqref{frame.operators.spaces} that $\mathcal{M}\mathcal{A}=S_{X}$. \qed
\begin{theorem}\cite[corollary 1, p.248]{pel2}.
\label{uni.basis.1969.C1}
There exists a separable Banach space $\mathcal{B}$ (resp. $\tilde{\mathcal{B}}$) with 
a Schauder  basis (an unconditional Schauder basis) such that every separable Banach 
space with a Schauder basis (an unconditional Schauder basis) is isomorphic to a 
complemented subspace of $\mathcal{B}$ (resp. of $\tilde{\mathcal{B}}$).
\end{theorem}
\begin{corollary}
\label{universal.1}
There exists a  Banach space $\mathcal{B}$  with 
a Schauder basis  such that every Banach 
space with a Schauder frame is isomorphic to a complemented subspace of $\mathcal{B}$.
\end{corollary}
\textbf{Proof.}
Since  $S_{X}$ is the image of $\mathcal{A}$ by the bounded projection 
$\mathcal{M}$ then $S_{X}$ is complemented in $\mathcal{A}$ 
\cite[corollary 3.2.14, p.299]{meg}.  According to theorem \eqref{uni.basis.1969.C1}
there exists a  Banach space $\mathcal{B}$
with a Schauder basis such that $\mathcal{A}$ 
is isomorphic to a  complemented subspace of $\mathcal{B}$, and  since $X$ is
isomorphic to $S_{X}$ then, $X$ is isomorphic to a complemented subspace of $\mathcal{B}$    \qed
\begin{proposition} 
\label{coro.complemented.frame.1}
Let $X$  be a Banach space with a Schauder frame.  If
$X_{0}$ a    complemented subspace in  $X$
 then  $X_{0}$   has a Schauder frame.
\end{proposition}
\textbf{Proof.}
Assume that $X_{0}$ is a complemented subspace in $X$, then there exists a continuous 
projection $P:X\longrightarrow X_{0}$.
 Let $\mathcal{F}$ be a Schauder frame of $X$, then we have for each 
$x\in X_{0}$:
\begin{align*}
x =&P\left( x\right) \\
=&\overset{+\infty}{\underset{n=1}{\sum }}
f_{n}(x)P\left(x_{n}\right)
\end{align*}
It follows that the paire 
$\left\lbrace \left( P\left(x_{n}\right),
f_{n\left \vert X_{0}\right. } \right) \right\rbrace _{n\in \mathbb{N}^{*}}$
 is a Schauder frame of $X_{0}$ (where $f_{n\left \vert X_{0}\right. } $
 is the restriction of $f_{n} $ to $X_{0}$ for each $n\in \mathbb{N}^{*}$).\qed
\begin{corollary}
\label{characterization.conditional.1}
 There exists a separable Banach space $\mathcal{B}$ 
with a Schauder basis  such that a Banach space  $X$ 
has a Schauder frame if and only if  $X$  is isomorphic to a
complemented subspace  of  $\mathcal{B}$.
\end{corollary}
\textbf{Proof.}
Il follows directly from corollary \eqref{universal.1}
and   proposition \eqref{coro.complemented.frame.1}.\qed
\begin{lemma}\cite[theorem 3.13, p.290]{W.B.Johnson.J.Lindenstrauss.2001}
\label{Characterization.BAP.Complemented}
A separable Banach space $X$ has the bounded approximation property
 if and only if
it is isomorphic to a complemented subspace of a Banach space with
 Schauder basis.
\end{lemma}
\begin{theorem}
\label{Schauder frame.eq.BAP}
A separable Banach space $X$ has a Schauder frame if and only if
it has the bounded approximation property. 
\end{theorem}
\textbf{Proof.}
Il follows directly from corollary \eqref{characterization.conditional.1}
and   lemma \eqref{Characterization.BAP.Complemented}.\qed
\begin{example}
Let $(1\leq p \leq \infty )$. The Banach space 
$\left (\sum_{n=1}^{\infty}L_{p}[0,1]\right )_{p}$ 
 has a Schauder frame.
\end{example}
 \textbf{Proof.}
The Banach space $L_{p}[0,1]$ has a Schauder basis, then has a Schauder frame, 
and consequently by theorem \eqref{Schauder frame.eq.BAP}, 
it has the bounded approximation property. 
 Il follows  from \cite[p.287]{W.B.Johnson.J.Lindenstrauss.2001}
 that $\left (\sum_{n=1}^{\infty}L_{p}[0,1]\right )_{p}$ has the 
 bounded approximation property. Consequently, also by
 theorem \eqref{Schauder frame.eq.BAP}, 
 $\left (\sum_{n=1}^{\infty}L_{p}[0,1]\right )_{p}$ 
 has a Schauder frame. \qed
 
The result of proposition \eqref{tensor.frames} already proved
 in \cite{kabbaj.karkri.zoubeir.2023}  by means of another way.
  Now, we give a simple proof of this result.
 \begin{proposition}
 \label{tensor.frames}
If $X$ and $Y$ are Banach spaces  with Schauder frames then, the
 Banach space $  X\widehat{\otimes}_{\pi}Y$
 (the projective tensor product of $X$ and $Y$ \cite{Raymond.A.Ryan}) has a Schauder frame.
\end{proposition}
\textbf{Proof.}
Since $X$ and $Y$ have Schauder frames, 
then by theorem \eqref{Schauder frame.eq.BAP},
$X$ and $Y$ have the bounded approximation property.  
 By means of  \cite[exercice 4.5, p.92]{Raymond.A.Ryan}
 $  X\widehat{\otimes}_{\pi}Y$ has the bounded approximation property.
 Finally, also by theorem \eqref{Schauder frame.eq.BAP}, 
 $X$ and $Y$ have Schauder frames.\qed
\begin{proposition}
There exists a separable Banach space without Schauder basis
 which has a Schauder frame.
\end{proposition}
\textbf{Proof.}
Thanks to S.J.Szarek \cite{S.J.Szarek} there exist 
reflexive Banach spaces $X$, $Y$ such that both
 $Y$ and $X\oplus Y$ have Schauder bases  but $X$ does not. Then 
 $X\oplus Y$ has a Schauder frame. On the other hand, since $X$
is complemented in $X\oplus Y$
 it  has a Schauder frame (proposition \eqref{coro.complemented.frame.1}).\qed
\begin{corollary}
There exists a separable Banach space without Schauder frame.
\end{corollary}
\textbf{Proof.}
Thank to Enflo  \cite{P.Enflo.1973}, there exists a separable Banach space
$X_{E}$ which fails to have the approximation property, 
then $X_{E}$ fails to have the bounded approximation property.  Consequently, by   
theorem \eqref{Schauder frame.eq.BAP}, 
$X_{E}$ has no Schauder frame.\qed
\begin{example}
The Banach space $\mathcal{L}(\mathcal{H},\mathcal{H})$ of all 
bounded linear operators on a Hilbert space $\mathcal{H}$  has no Schauder frame.
\end{example}
\textbf{Proof.}
Thanks to \cite{A.Szankowski}
the Banach space $\mathcal{L}(\mathcal{H},\mathcal{H})$
 has no the approximation property, then has no the bounded approximation property.
 Consequently, by 
 theorem \eqref{Schauder frame.eq.BAP}, $\mathcal{L}(\mathcal{H},\mathcal{H})$
  has no Schauder frame.\qed
\section{Unconditional and besselian Schauder frames}
\begin{lemma}\cite[lemma 1.45, p.31]{A.Abramovich.2002}
\label{equi.unconditional.cv}
For a series $\sum x_{n}$ in Banach space $X$ the following are 
equivalent.
\begin{enumerate}
\item The series $\sum x_{n}$ converges unconditionally.
\item
For any sequence 
$  \left( s_{n}\right)_{n\in \mathbb{N}^{*}}$  of signs
(i.e., $s_{n}=\pm 1$ for each $n$) the series
$\sum s_{n}x_{n}$  is norm convergent.
\item
For each $\varepsilon >0$ there exists a natural number  $k$
such that for each finite subset $A$ of $\mathbb{N}$ with $min(A)\geq k$ we have
$\left \Vert \sum_{n\in A}x_{n}  \right \Vert_{X} < \varepsilon$.
\end{enumerate}
\end{lemma}
\begin{lemma}
Assume that $\mathcal{F}$ is an unconditionally Schauder frame of $X$.
Given an arbitrary sign sequence 
$  s=\left\lbrace s_{n}\right\rbrace _{n\in \mathbb{N}^{*}}\in \mathcal{S}$. 
The mapping $M_{s}:X\longrightarrow X$ given by the formula
\[M_{s}(x)=\underset{j=1}{\overset{+\infty }{\sum }}
s_{j} f_{j}(x)x_{j} \]
is a linear bounded operator and we have
\[k:=\underset{s\in \mathcal{S}}{\sup}\left \Vert  M_{s} \right \Vert<\infty\]
\end{lemma}
\textbf{Proof.} Is similar to that used in the proof in
\cite[p.32-33]{A.Abramovich.2002} for an unconditionally Schauder basis.\qed
\begin{proposition}
\label{111}
Assume that $\mathcal{F}$ is an unconditionally Schauder frame of $X$.
Then $\mathcal{F}$ is a besselian Schauder frame of $X$.
\end{proposition}
\textbf{Proof.}
Let $x\in X$ and $f\in \mathbb{B}_{X^{*}}$. If 
$  \left\lbrace s_{n}\right\rbrace _{n\in \mathbb{N}^{*}}$ is a sequence of signs
such that $ s_{n}=sign\left(f_{n}(x)f(x_{n}) \right)$ for each
$n\in \mathbb{N}^{*}$ then
\begin{align*}
\underset{j=1}{\overset{+\infty }{\sum }}
\left \vert f_{j}(x)\right\vert \left \vert f(x_{j})\right \vert 
& =\overset{+\infty }{\underset{j=1}{\sum }}
s_{j}f_{j}(x) f\left( x_{j}\right)\\
&=\left \vert f\left(\overset{+\infty }{
\underset{j=1}{\sum }}
s_{j}f_{j}(x) x_{j}\right)\right \vert\\
&=\left \vert f\left(M_{s}(x)\right)\right \vert\\
&\leq \left \Vert  M_{s}(x) \right \Vert\\
&\leq \left \Vert  M_{s} \right \Vert \left \Vert  x \right \Vert\\
&\leq k \left \Vert  x \right \Vert
\end{align*}
Consequently, we have
\[\underset{j=1}{\overset{+\infty }{\sum }}
\left \vert f_{j}(x)\right
\vert \left \vert f(x_{j})\right \vert 
\leq k\left \Vert  x \right \Vert_{X} \left \Vert  f \right \Vert_{X^{*}} 
\;\;\; (x\in X, f\in X^{*})\]
\qed
\begin{corollary}
\label{sup.exists}
Assume that $\mathcal{F}$ is an unconditionally Schauder frame of $X$. Then
for each $x\in X$, we have
\[\underset{\sigma \in \mathfrak{S}_{\mathbb{N}^{* }}, n\in \mathbb{N}^{*}}{\sup }
\left \Vert \overset{n}{\underset{j=1}{\sum }}
f_{\sigma (j)}(x)x_{\sigma (j)}\right \Vert _{X}
\leq \mathcal{L_{\mathcal{F}}}\Vert x\Vert _{X}\] 
\end{corollary}
\textbf{Proof.}
\begin{align*}
\underset{\sigma \in \mathfrak{S}_{\mathbb{N}^{*}}, n\in \mathbb{N}^{*}}{\sup }
\left \Vert \overset{n}{\underset{j=1}{\sum }}
f_{\sigma (j)}(x)x_{\sigma (j)}\right \Vert _{X} 
& =\underset{\sigma \in \mathfrak{S}_{\mathbb{N}^{*}}, n\in \mathbb{N}^{*}}{\sup }
\left( \underset{f\in \mathbb{B}_{X^{*}}}{\sup }
\left \vert \overset{n}{\underset{j=1}{\sum }}
f_{\sigma (j)}(x)f(x_{\sigma (j)})\right \vert \right) \\
& \leq \underset{\sigma \in \mathfrak{S}_{\mathbb{N}^{*}}}{\sup }
\left(\underset{f\in \mathbb{B}_{X^{*}}}{\sup }\overset{+\infty }{\underset{n=1}{\sum }}
\left \vert f_{\sigma (n)}(x)f(x_{\sigma (n)})\right \vert \right) \\
& \leq \underset{f\in \mathbb{B}_{X^{*}}}{\sup }
\overset{+\infty }{\underset{n=1}{\sum }}
\left \vert f_{n}(x)f(x_{n})\right
\vert \\
& \leq \mathcal{L_{\mathcal{F}}}\Vert x\Vert _{X}
\end{align*}
\qed
\begin{proposition}\cite{karkri.zoubeir.2021}.
\label{2}
Assume that  $X$  is a weakly
sequentially complete Banach space and that  $\mathcal{F}$ 
is a besselian Schauder frame of  $X $. Then for each $x\in
X$, the series $\sum_{n} f_{n}\left( x\right) x_{n}$ 
is unconditionally convergent  to  $x$.
\end{proposition} 
\begin{corollary}
Assume that  $X$  is a weakly
sequentially complete Banach space and  
$\mathcal{F}$ 
be a paire of  $X$. Then $\mathcal{F}$ is unconditionally
Schauder frame of $X$ if and only if $\mathcal{F}$ 
is a besselian Schauder frame of $X$.
\end{corollary}
\begin{proposition}
Let $\left\lbrace x_{n} \right\rbrace_{n\in\mathbb{N}^{* }}$
be a sequence in a Banach space $X$, 
such that $x_{n}\neq 0$ $(n=1,2,...)$ and let 
$\tilde{\mathcal{A}}$ be the linear space
of sequences of scalars
\[\tilde{\mathcal{A}}=\left\lbrace
 \left\lbrace \alpha_{n} \right\rbrace_{n\in\mathbb{N}^{*}}\subset \mathbb{K}:
\sum _{k}\alpha_{k}x_{k}\;converges\; unconditionally \right\rbrace\]
endowed with the norm
\begin{align}
\label{norm.unconditional}
\left \Vert \left\lbrace \alpha_{n} \right\rbrace_{n\in\mathbb{N}^{* }}
 \right \Vert_{\tilde{\mathcal{A}}}
=\underset{\sigma \in \mathfrak{S}_{\mathbb{N}^{*}}}{\sup }
\underset{n\in\mathbb{N}^{*}}{\sup }
\left \Vert \overset{n}{\underset{k=1}{\sum }}
 \alpha_{\sigma(k)}x_{\sigma(k)}  \right \Vert_{X}
\end{align}
Then $\tilde{\mathcal{A}}$ is a Banach space.
\end{proposition}

\textbf{Proof.}
Let $\varepsilon >0$. Since the series $\sum_{j} \alpha_{j}x_{j}$ is unconditionally
convergent in $X$, there exists a positive
integer $k(\varepsilon )\in \mathbb{N}^{*}$ such that the inegality 
\begin{equation*}
\left \Vert \sum_{j\in L}\alpha_{j}x_{j}\right \Vert _{X}\leq
\varepsilon
\end{equation*}
holds for each $L\in \mathcal{D}_{\mathbb{N}^{*}}$ with
 $\min \left( L\right) \geq k(\varepsilon )$.
 Then for each  $L\in \mathcal{D}_{\mathbb{N}^{*}}$ we have
 $$\left \Vert \underset{j\in L\cap \left \{ 1,...k(\varepsilon)\right \} }{\sum }
\alpha_{j}x_{j}\right \Vert _{X}
\leq \overset{k(\varepsilon)}{\underset{j=1}{\sum }}
\left \Vert \alpha_{j}x_{j}\right \Vert _{X}
\;\;and \;\; \left \Vert \underset{j\in L\backslash \left \{ 1,...k(\varepsilon)\right \} }{\sum }
\alpha_{j}x_{j}\right \Vert _{X}\leq \varepsilon$$
It follows that: 
\begin{equation*}
\left \Vert \underset{j\in L}{\sum }\alpha_{j}x_{j}\right \Vert_{X}
\leq 
\varepsilon+\overset{k(\varepsilon)}{\underset{j=1}{\sum }}
\left \Vert\alpha_{j} x_{j}\right \Vert _{X}
\end{equation*}
Consequently, 
\begin{equation*}
\underset{L\in \mathcal{D}_{\mathbb{N}^{*}}}{\sup }
\left \Vert\sum_{j\in L} \alpha_{j} x_{j}\right \Vert _{X}<+\infty
\end{equation*}
then 
\begin{equation*}
\underset{\sigma \in \mathfrak{S}_{\mathbb{N}^{*}},\;n\in \mathbb{N}^{*}}{\sup }
\left \Vert \underset{j=1}{\overset{n}{\sum }}
\alpha_{\sigma (j)}x_{\sigma (j)}\right \Vert_{X}<+\infty
\end{equation*}
Hence the quantity \eqref{norm.unconditional} is finite. \\
Let $\left\lbrace 
\left\lbrace \alpha_{n}^{k}\right\rbrace_{n\in\mathbb{N}^{*}}
 \right\rbrace_{k\in\mathbb{N}^{*}} $
  be a Cauchy sequence in $\tilde{\mathcal{A}}$.
 Then for each $\varepsilon>0$ there exists a positive integer $N(\varepsilon)$
 such that
\[\left\Vert \left\lbrace \alpha_{n}^{k}\right\rbrace_{n\in\mathbb{N}^{*}}
-\left\lbrace \alpha_{n}^{m}\right\rbrace_{n\in\mathbb{N}^{*}}\right \Vert_{\tilde{\mathcal{A}}}
= \underset{\sigma \in \mathfrak{S}_{\mathbb{N}^{*}},\;n\in \mathbb{N}^{*}}{\sup }
\left \Vert \underset{j=1}{\overset{n}{\sum }}
\left (\alpha_{\sigma (j)}^{k}-\alpha_{\sigma (j)}^{m}\right )x_{\sigma (j)}\right \Vert_{X}<\varepsilon
\;\; (k,m>N(\varepsilon))\]
hence, for each $k,m>N(\varepsilon)$, $n\in\mathbb{N}^{*}$
 and $\sigma\in \mathfrak{S}_{\mathbb{N}^{*}}$ 
 we have
$$\left\Vert \left ( \alpha_{\sigma (n)}^{k}-\alpha_{\sigma (n)}^{m}\right )x_{\sigma (n)}\right\Vert_{X}
\leq \left \Vert \underset{j=1}{\overset{n}{\sum }}
\left ( \alpha_{\sigma (j)}^{k}-\alpha_{\sigma (j)}^{m}\right )x_{\sigma (j)}\right \Vert_{X}
+
\left \Vert \underset{j=1}{\overset{n-1}{\sum }}
\left ( \alpha_{\sigma (j)}^{k}-\alpha_{\sigma (j)}^{m}\right )x_{\sigma (j)}\right \Vert_{X}
\leq 2\varepsilon$$
it follows that 
\[\left \vert \alpha_{\sigma (n)}^{k}-\alpha_{\sigma (n)}^{m}\right \vert
<\frac{2\varepsilon}{\left \Vert x_{\sigma (n)} \right \Vert_{X}}\;\;
\left (k,m>N(\varepsilon), n\in\mathbb{N}^{*}
 \;and \;\sigma\in \mathfrak{S}_{\mathbb{N}^{*}}\right )\]
Consequently, for each $n\in\mathbb{N}^{*}$ and
 $\sigma\in \mathfrak{S}_{\mathbb{N}^{*}}$ the sequence
$\left\lbrace \alpha_{\sigma (n)}^{k}\right\rbrace_{k\in\mathbb{N}^{*}}$
convergent to a scalar $\alpha_{\sigma (n)}$.
 From the inequalities
\[\left \Vert \underset{j=1}{\overset{n}{\sum }}
\left (\alpha_{\sigma (j)}^{k}-\alpha_{\sigma (j)}^{m}\right )x_{\sigma (j)}\right \Vert_{X}
<\varepsilon
\;\; \left (k,m>N(\varepsilon), n\in\mathbb{N}^{*},
 \sigma\in \mathfrak{S}_{\mathbb{N}^{*}}\right ) \]
 we obtain, 
 for $m\longrightarrow \infty$,
 \[\left \Vert \underset{j=1}{\overset{n}{\sum }}
\left (\alpha_{\sigma (j)}^{k}-\alpha_{\sigma (j)}\right )x_{\sigma (j)}\right \Vert_{X}
<\varepsilon
\;\; \left (k>N(\varepsilon), n\in\mathbb{N}^{*},
 \sigma\in \mathfrak{S}_{\mathbb{N}^{*}}\right ) \]
Then, for each 
$k>N(\varepsilon)$, $n,r\in\mathbb{N}^{*}$ and
 $\sigma\in \mathfrak{S}_{\mathbb{N}^{*}}$ we have 
\begin{align*}
\left \Vert \underset{j=n+1}{\overset{n+r}{\sum }}
 \alpha_{\sigma (j)}x_{\sigma (j)}\right \Vert_{X}
&\leq \left \Vert \underset{j=n+1}{\overset{n+r}{\sum }}
\left (\alpha_{\sigma (j)}^{k}-\alpha_{\sigma (j)}\right )x_{\sigma (j)}\right \Vert_{X}
+\left \Vert \underset{j=n+1}{\overset{n+r}{\sum }}
 \alpha_{\sigma (j)}^{k}x_{\sigma (j)}\right \Vert_{X}\\
 \leq & 2\varepsilon+\left \Vert \underset{j=n+1}{\overset{n+r}{\sum }}
 \alpha_{\sigma (j)}^{k}x_{\sigma (j)}\right \Vert_{X}
\end{align*}
whence, since each series $\sum_{j}\alpha_{j}^{k}x_{j}$ 
is unconditionally convergent, and since $X$ is complete, it follows that
 $\sum_{j}\alpha_{j}x_{j}$ converges unconditionally, that is 
 $\left\lbrace \alpha_{n} \right\rbrace_{n\in\mathbb{N}^{*}} \in\tilde{A}$. 
Moreover, by the above we have 
\[\left\Vert \left\lbrace \alpha_{n}^{k}\right\rbrace_{n\in\mathbb{N}^{*}}
-\left\lbrace \alpha_{n}\right\rbrace_{n\in\mathbb{N}^{*}}\right \Vert_{\tilde{\mathcal{A}}}
= \underset{\sigma \in \mathfrak{S}_{\mathbb{N}^{*}},\;n\in \mathbb{N}^{*}}{\sup }
\left \Vert \underset{j=1}{\overset{n}{\sum }}
\left (\alpha_{\sigma (j)}^{k}-\alpha_{\sigma (j)}\right )a_{\sigma (j)}\right \Vert_{X}<\varepsilon
\;\; (k>N(\varepsilon))\]
which completes the proof of the proposition.\qed
\begin{remark}
\label{frame.uncond.SX}
Assume that $\mathcal{F}$ is an unconditionally Schauder frame of $X$. Then
the Banach space $S_{X}$ is a subspace of $\tilde{\mathcal{A}}$. 
\end{remark}
\textbf{Proof.}
It follows directly from corollary \eqref{sup.exists}.\qed
\begin{remark}
If  $\left\lbrace x_{n} \right\rbrace_{n\in\mathbb{N}^{*}}$
is a sequence in a Banach space $X$ such that $x_{n}\neq 0$ $(n=1,2,...)$, 
then the mapping 
\begin{align*}
\tilde{T} :\tilde{\mathcal{A}} &\rightarrow  X \\ 
 {\left\lbrace \alpha_{n} \right\rbrace}  &\mapsto \overset{+\infty }{\underset{n=1}{\sum }}
 \alpha_{k}x_{k}
\end{align*}
defines a bounded linear operator.
\end{remark}
\textbf{Proof.}
For each $\alpha=\left\lbrace \alpha_{n} \right\rbrace_{n\in\mathbb{N}^{*}}
\in \mathbb{B}_{\mathcal{A}}$ we have
 $$\left \Vert \sum_{k=1}^{n}\alpha_{\sigma(k)}x_{\sigma(k)} \right \Vert_{X}
 \leq 1 \;\;\;(n\in\mathbb{N}^{*}, \sigma\in \mathfrak{S}_{\mathbb{N}^{*}})$$
 Consequently,
\begin{align*}
\left \Vert \tilde{T} \right \Vert
&=\underset{\alpha\in\mathbb{B}_{\mathcal{A}}}{\sup }
\left \Vert \tilde{T}\alpha \right \Vert_{X}\\
 &=\underset{\alpha\in\mathbb{B}_{\mathcal{A}}}{\sup }
\left \Vert \overset{\infty}{\underset{k=1}{\sum }} \alpha_{k}x_{k}  \right \Vert_{X}\\
&\leq 1
\end{align*}
\qed
\begin{remark}
Let  $\mathcal{F}$ 
be a Schauder frame of a Banach space $X $. For each $k\in\mathbb{N}^{* }$
and $ \sigma \in\mathfrak{S}_{\mathbb{N}^{*}}$  the  mapping 
\begin{align*}
T_{k,\sigma}:\tilde{\mathcal{A}} &\rightarrow  X \\ 
 \alpha={\left\lbrace \alpha_{n} \right\rbrace} 
  &\mapsto \overset{k}{\underset{i=1}{\sum }}
 f_{\sigma(i)}(\tilde{T}\alpha)x_{\sigma(i)}
\end{align*}
defines a bounded linear operator.
\end{remark}
\textbf{Proof.}
For each $k\in\mathbb{N}^{*}$ and
$ \sigma \in\mathfrak{S}_{\mathbb{N}^{*}}$, $T_{k,\sigma}$ 
is a finite sum of bounded linear operators
$$\alpha \longmapsto f_{\sigma(i)}(\tilde{T}\alpha)x_{\sigma(i)}, \;\;\;(1\leq i\leq k)$$\qed

\begin{proposition}
The unit vectors 
$e_{n}=\left\lbrace \delta_{n,k} \right\rbrace_{k\in\mathbb{N}^{* }}$ 
$(n=1,2,...)$ constitue an unconditional Schauder basis of  
$\tilde{\mathcal{A}}$.
\end{proposition}
\textbf{Proof.}
Let $\left\lbrace \alpha_{i} \right\rbrace_{i\in\mathbb{N}^{* }}
\in \tilde{\mathcal{A}}  $ 
 and $\sigma _{0}\in \mathfrak{S}_{\mathbb{N}^{*}}$. 
 We have by \eqref{norm.unconditional}:
\begin{align*}
 \left \Vert \left\lbrace \alpha_{i} \right\rbrace_{i\in\mathbb{N}^{*}}
 -\overset{n}{\underset{k=1}{\sum }}\alpha_{\sigma _{0}(k)}e_{\sigma _{0}(k)}
 \right \Vert_{\tilde{\mathcal{A}}}
 =&\underset{\underset{\sigma \in \mathfrak{S}_{\mathbb{N}^{*}}}{r\in \mathbb{N}^{* }}}{\sup }
 \left \Vert \underset{\underset{\sigma (j)\in 
\mathbb{N}^{*}\backslash \left \{ \sigma _{0}(1),...,\sigma
_{0}(n)\right \}}{1\leq j\leq r  }}{\sum }
\alpha_{\sigma(j)}x_{\sigma(j)}\right \Vert_{X} \;\;\;(n\in\mathbb{N}^{*})
\end{align*}
Let $\varepsilon >0$. Since the series 
$\sum_{j}\alpha_{j}x_{j}$ is unconditionally convergent in $X$, 
by lemma \eqref{equi.unconditional.cv} 
 there exists an integer $k(\varepsilon )\in \mathbb{N}^{*}$ such that: 
\begin{equation*}
\left \Vert \underset{j\in L}{\sum }\alpha_{j}x_{j}\right \Vert_{X}
< \varepsilon
\end{equation*}
for each $L\in \mathcal{D}_{\mathbb{N}^{* }}$ with $\min (L)\geq
k(\varepsilon )$. Let $n\in \mathbb{N}^{* }$ such that
 $$ \left\lbrace 1,2,...,k(\varepsilon ) \right\rbrace 
\subset \left\lbrace \sigma_{0}(1),\sigma_{0}(2),...,\sigma_{0}(n) \right\rbrace$$
then  $\mathbb{N}^{* }\backslash \left\lbrace \sigma_{0}(1),\sigma_{0}(2),...,\sigma_{0}(n) \right\rbrace
\subset \mathbb{N}^{* }\backslash \left\lbrace 1,2,...,k(\varepsilon ) \right\rbrace$
and consequently
\[\left \Vert \underset{\underset{\sigma (j)\in 
\mathbb{N}^{*}\backslash \left \{ \sigma _{0}(1),...,\sigma
_{0}(n)\right \}}{1\leq j\leq r  }}{\sum }
\alpha_{\sigma(j)}x_{\sigma(j)}\right \Vert_{X}  
< \varepsilon\;\;\; (\sigma\in \mathfrak{S}_{\mathbb{N}^{*}},  r\in \mathbb{N}^{*})
\]
then 
\begin{equation*}
\underset{n\rightarrow +\infty }{\lim }
 \left \Vert \left\lbrace \alpha_{i} \right\rbrace_{i\in\mathbb{N}^{*}}
 -\overset{n}{\underset{k=1}{\sum }}\alpha_{\sigma _{0}(k)}e_{\sigma _{0}(k)}
 \right \Vert_{\tilde{\mathcal{A}}}=0
\end{equation*}
\qed
\begin{theorem}
Let $\mathcal{F}$ 
be an unconditional Schauder frame of a Banach space $X $. Then 
the infinite matrix
\[\tilde{\mathcal{M}}_{i,j}=f_{i}\left( x_{j}\right ),\; (i,j\in \mathbb{N}^{* })\]
defines a bounded projection $\tilde{\mathcal{M}}$ from $\mathcal{\tilde{A}}$ onto $S_{X}$. 
\end{theorem}
\textbf{Proof.}
Let 
$ \alpha=\left\lbrace \alpha_{n} \right\rbrace_{n\in\mathbb{N}^{*}}$
be a sequence in $\mathcal{\tilde{A}}$. With the same argument as 
that used in the proof of theorem \eqref{projection.onto} and  
the fact that the series $\sum_{k}f_{k}\left(\tilde{T}\alpha \right )x_{k}$ 
converge unconditionally to $\tilde{T}\alpha$, we obtain 
$$\tilde{\mathcal{M}}\alpha=\left\lbrace f_{n}\left(\tilde{T}\alpha \right )
 \right\rbrace_{n\in\mathbb{N}^{*}}\in S_{X}\subset \tilde{\mathcal{A}} $$
also $\tilde{\mathcal{M}}$ is a projection in $\tilde{\mathcal{A}}$. By  the
uniform boundedness principle applied  to the family of bounded linear operators
$\left\lbrace T_{n,\sigma} \right\rbrace_{n\in \mathbb{N}^{*}, 
\sigma \in\mathfrak{S}_{\mathbb{N}^{*}}}$
and  corollary \eqref{sup.exists} we have
 \begin{align*}
 \left \Vert \tilde{\mathcal{M}} \right \Vert
 &=\underset{\alpha\in\mathbb{B}_{\mathcal{\tilde{A}}}}{\sup }
 \left \Vert \tilde{\mathcal{M}}\alpha \right \Vert_{\mathcal{A}}\\
 &=\underset{\alpha\in\mathbb{B}_{\mathcal{A}}}{\sup }
 \left \Vert \left\lbrace f_{n}\left(\tilde{T}\alpha \right )
 \right\rbrace_{n\in\mathbb{N}^{*}} \right \Vert_{\mathcal{A}}\\
 &=\underset{\alpha\in\mathbb{B}_{\tilde{\mathcal{\tilde{A}}}}}{\sup }
 \underset{n\in\mathbb{N}^{*}}{\sup }
 \underset{\sigma \in\mathfrak{S}_{\mathbb{N}^{*}}}{\sup }
 \left \Vert \overset{n}{\underset{k=1}{\sum }}
f_{\sigma(k)}\left(\tilde{T}\alpha \right )x_{\sigma(k)} \right \Vert_{X}\\
&=\underset{\alpha\in\mathbb{B}_{\tilde{\mathcal{A}}}}{\sup}
 \underset{n\in\mathbb{N}^{*}}{\sup }
 \underset{\sigma \in\mathfrak{S}_{\mathbb{N}^{*}}}{\sup}
 \left \Vert T_{n,\sigma}\alpha  \right \Vert_{X}\\
 &=\underset{n\in\mathbb{N}^{*}}{\sup }
  \underset{\sigma \in\mathfrak{S}_{\mathbb{N}^{*}}}{\sup}
 \underset{\alpha\in\mathbb{B}_{\tilde{\mathcal{A}}}}{\sup}
 \left \Vert T_{n,\sigma}\alpha  \right \Vert_{X}\\
 &=\underset{n\in\mathbb{N}^{*}}{\sup}
 \underset{\sigma \in\mathfrak{S}_{\mathbb{N}^{*}}}{\sup}
 \left \Vert T_{n,\sigma}\right \Vert <\infty
 \end{align*}
then $\tilde{\mathcal{M}}$ is a bounded operator. Finally, for each $n\in  \mathbb{N}^{*} $ we have
$\tilde{\mathcal{M}}e_{n}=\left\lbrace f_{i}\left(x_{n} \right ) \right\rbrace_{i\in\mathbb{N}^{*}}$,
consequently, it follows from remark \eqref{frame.operators.spaces} 
that $\tilde{\mathcal{M}}\mathcal{\tilde{A}}=S_{X}$. \qed
\begin{corollary}
\label{universal.2}
There exists a separable Banach space $\tilde{\mathcal{B}}$  with 
an unconditional basis  such that every separable Banach 
space with an unconditional  Schauder frame is isomorphic 
to a complemented subspace of $\tilde{\mathcal{B}}$.
\end{corollary}
\textbf{Proof.}
Since  $S_{X}$ is the image of $\tilde{ \mathcal{A}}$ 
by the bounded projection 
$\tilde{\mathcal{M}} $ then $S_{X}$ is complemented in $\tilde{ \mathcal{A}}$ 
\cite[corollary.3.2.14, p.299]{meg}. According to theorem \eqref{uni.basis.1969.C1}
there exists a separable Banach space $\tilde{\mathcal{B}}$
with an unconditional Schauder basis such that $\tilde{ \mathcal{A}}$ 
is isomorphic to a  complemented subspace of $\tilde{\mathcal{B}}$, and  since $X$ is
isomorphic to $S_{X}$ then, $X$ is isomorphic to a complemented subspace 
of $\mathcal{\tilde{\mathcal{B}}}$.    \qed

It is an open problem whether every complemented subspace of a space with
unconditional Schauder basis has an unconditional Schauder basis
 \cite[p.279]{W.B.Johnson.J.Lindenstrauss.2001}.
In particular, it is unknown whether every separable $\mathcal{L}_{p}$,
$1<p<\infty$ $(p\neq 2)$ space has an unconditional basis
 \cite[p.59]{W.B.Johnson.J.Lindenstrauss.2001}.
 The proposition \eqref{coro.complemented.frame.2} 
 give a positive answer to this problem
  (for the setting of unconditional Schauder frames),
 that  is,  an unconditional Schauder frame passes to complemented subspaces.
\begin{proposition} 
\label{coro.complemented.frame.2}
Let $X$  be a Banach space with an unconditional Schauder frame   and
$X_{0}$ a subspace of $X$  which is complemented in  $X$.
 Then  $X_{0}$   has an unconditional Schauder frame.
\end{proposition}
\textbf{Proof.}
Assume that $X_{0}$ is a complemented subspace in $X$. Then there exists a continuous 
projection $P:X\longrightarrow X_{0}$.
 Let $\mathcal{F}$ be an unconditional Schauder frame of $X$. 
Then we have for each 
$x\in X_{0}$:
\begin{align*}
x =&P\left( x\right) \\
=&\overset{+\infty}{\underset{n=1}{\sum }}
f_{\sigma(n)}(x)P\left(x_{\sigma(n)}\right),\;\;\;(\sigma \in \mathfrak{S}_{\mathbb{N}^{* }})
\end{align*}
It follows that the paire 
$\left\lbrace \left( P\left(x_{n}\right),
f_{n\left \vert X_{0}\right. }\right) \right\rbrace _{n\in \mathbb{N}^{*}}$
 is an unconditional Schauder frame of $X_{0}$ 
 where $f_{n\left \vert X_{0}\right. }$
 is the restriction of $f_{n}$ to $X_{0}$ for each $n\in \mathbb{N}^{*}$.\qed
\begin{example}
If $1<p<\infty$ and $X$ is a separable $\mathcal{L}_{p}$ space, 
then $X$ has an unconditional Schauder frame.
\end{example}
\textbf{Proof.}
According to \cite[p.57]{W.B.Johnson.J.Lindenstrauss.2001}, 
$X$ is  isomorphic to a complemented subspace of $L_{p}[0,1]$
 which has an unconditional Schauder frame, then
    by proposition \eqref{coro.complemented.frame.2}
   $X$ has an unconditional Schauder frame.\qed
\begin{corollary}
There exists a separable Banach space  $\tilde{\mathcal{B}}$
   with an unconditional  Schauder basis
      such that a Banach space  $X$ 
      has an unconditional Schauder frame if and only if  $X$
  is isomorphic to a complemented subspace of  $\tilde{\mathcal{B}}$.
\end{corollary}
\textbf{Proof.}
Il follows directly from corollary \eqref{universal.2}
and   proposition \eqref{coro.complemented.frame.2}.\qed

\bigskip

\bigskip

\end{document}